\documentclass[12pt]{article}
\usepackage[latin1]{inputenc}
\usepackage{amsfonts}
\usepackage{amsmath}
\usepackage{amssymb,amscd,latexsym}
\usepackage{enumerate}
\usepackage[dvips]{graphicx}
\newtheorem{lemma}{Lemma}
\newtheorem{theorem}{Theorem}
\newtheorem{prop}{Proposition}
\newtheorem{cor}{Corollary}

\newtheorem{definition}{Definition}

\newcommand{\implica}{\Rightarrow}

\newcommand{\erre}{\mathbb{R}}
\newcommand{\dem}{\noindent{\bf Proof: }}
\newcommand{\sen}{\mathrm{sen}}
\newcommand{\caixapreta}{\rule{2mm}{2mm}}
\newcommand{\cqd}{\begin{flushright}\caixapreta\vspace{5mm}\end{flushright}}

\newcommand{\te}{\theta}
\newcommand{\tet}[1]{\theta_{#1}}

\newcommand{\f}[1]{f(\theta_{#1})}
\newcommand{\ds}{\displaystyle}

\title{\vspace{-0.in}\parbox{\linewidth }{\footnotesize\noindent
} \\  \bf On stacked central configurations of the planar coorbital satellites problem}

\author{
\\ Allyson Oliveira \\
\small {Departamento de Matem\'atica, Universidade Federal de Pernambuco }
\\ \small{ Recife-PE,  CEP. 50540-740, Brazil. e-mail: allyson.oliveira@ufpe.br }
\\ \small{and}
\\ \small{N\'ucleo de Forma\c c\~ao Docente, Universidade Federal de Pernambuco}
\\ \small{Caruaru-PE, CEP. 55002-970, Brazil}
\\~\\ Hildeberto Cabral\\
\small {Departamento de Matem\'atica, Universidade Federal de Pernambuco }
\\ \small{ Recife-PE,  CEP. 50540-740, Brazil. e-mail: hild@dmat.ufpe.br}
}

\date{ }

\begin{document}
\maketitle

\author{ \ }

\begin{abstract}
In this work we look for central configurations of the planar $1+n$ body problem such that, after the addition of one or two satellites, we have a new planar central configuration. We determine all such configurations in two cases: the first, the addition of two satellites considering that all satellites have equal infinitesimal masses and the second case where one satellite is added but the infinitesimal masses are not necessarily equal.
\end{abstract}

\date{ }

\maketitle
\section{Introduction}

\hspace{0,5cm}
Central configurations of a system of $N$ bodies is one of the most classical and relevant topics in Celestial Mechanics. They are configurations such that the total Newtonian acceleration of every body is equal to a constant
multiplied by the position vector of this body with respect to the center
of mass of the configuration. Central configurations allow homographic motion, i. e. motion where the configuration of the system changes its size but keeps its shape. One of the reasons why central configurations are interesting is that they allow to obtain explicit homographic solutions of the N-body problem. They also arise as the limiting configuration of a total collapse.

The central configurations are invariant under rotation and homothety. So we are interested in the classes of central configurations modulus such transformations. There is an extensive literature concerning these configurations; see e.g. \cite{hagihara},\cite{saari}, \cite{wintner} for classic background.

In this work we study central configurations of the planar $1+n$ body problem, where we have one dominant mass and $n$ infinitesimal masses, called satellites, on a plane. This problem was first considered by Maxwell \cite{maxwell} trying to construct a model for Saturn's rings. Considering satellites with small and equal masses, Hall \cite{hall} shows that in a central configuration of $n$ satellites, the large body is at the center of a circle which passes through the satellites. Moreover, in this unpublished work, Hall shows that, if $n\geq e^{27000}$, there is a unique central configuration of this problem, that where the satellites are at the vertices of a regular polygon. Casasayas, Llibre and Nunes \cite{llibre1} proved
that the regular polygon is the only one if $n \geq e^{73}$. Cors, LLibre and Ollé \cite{llibre2} obtain numerically evidences that there is only one central configuration if $n\geq 9$ and that every central configuration is symmetric with respect to a straight line. Moreover they proved that there are only 3 symmetric central configurations of the $1+4$ body problem. Albouy and Fu \cite{albouy} proved that all central configurations of the $1+4$ body problem are symmetric which settles the question in this case.

In the case where the satellites do not have necessarily equal masses, Renner and Sicardy \cite{sicardy} obtained results about the inverse problem, that is, given a configuration of the coorbital satellites, find the infinitesimal masses making it a central configuration.
They also studied the linear stability. Corbera, Cors and LLibre \cite{llibre3} considering the $1+3$ body problem, found two different classes exhibiting symmetric and nonsymmetric configurations. And when two infinitesimal masses are equal, they provide evidence that the number of central configurations varies from five to seven.

Hampton \cite{hampton} provides a new family of planar central configurations
for the 5 body problem with the property that two bodies can
be removed and the remaining three bodies still form a central configuration. Such configurations are called stacked central configurations. Mello and Libre \cite{mello} studied one case of stacked planar central
configuration with 5 bodies in which three bodies are at the vertices of an equilateral triangle and the other two bodies lie on a mediatrix of the triangle.

In this work we study stacked central configurations of the planar $1+n$ body problem in two situations: the first, that of addition (or removal) of one satellite in the problem with arbitrary small masses and the case of addition (or removal) of two sattelites when all satellites have equal masses. Moreover the argument used in the latter case also shows that in a central configuration of the coorbital satellites problem if two satellites have the same mass the removal of them does not result in a central configuration if $n\geq 5.$

\section{Preliminaries}\label{preliminaries}

\hspace{0,5cm}Consider $N$ masses, $m_1,...,m_N$, in $\erre^2$ subject to their mutual Newtonian gravitational interaction. Let $M=diag\{m_1,m_1,...,m_N,m_N\}$ be the matrix of masses and let $q=(q_1,...,q_N), q_i\in \erre^2$ be the position vector. The equations of motion in an inertial reference frame with origin at the center of mass are given by
$$M\ddot{q} = \frac{\partial V}{\partial q},$$
where $V(q_1,...,q_N)=\ds\sum_{1\leq i<j\leq N}\frac{m_im_j}{\|q_i-q_j\|}.$

A non-collision configuration $q=(q_1,...,q_N)$ with $\sum_{i=1}^Nm_iq_i =0$ is a {\it central configuration} if there exists a positive constant $\lambda$ such that
$$M^{-1}V_q = \lambda q.$$

We are interested in the planar $N=1+n$ body problem, where the big mass is equal to 1 with position $q_0=0$. The remaining $n$ bodies with positions $q_i$, called satellites, have masses $m_i =\mu_i\epsilon, i=1,..,n$, where $\mu_i\in \erre^{+}$ and $\epsilon > 0$ is a small parameter that tends to zero.

In all central configuration of the planar $1+n$ body problem the satellites lie
on a circle centered at the big mass (\cite{llibre1,hall}), i.e. the satellites are coorbital. Since we are interested in central configuration modulus rotations and
homothetic transformations, we can assume that the circle has radius 1 and
that $q_1=(1,0)$.

We exclude collisions in the definition of central configuration and take as coordinates the angles $\theta_i$ between
two consecutive particles. See e.g. \cite{llibre1} for details. In this coordinates the space of configuration is the simplex $$\Delta = \{\theta = (\theta_1,...,\theta_n); \sum_{i=1}^n\theta_i = 2\pi, \theta_i>0, i=1,..,n\}$$ and the equations characterizing the central configurations of the planar 1+n body problem are
\begin{eqnarray}
&&\mu_2f(\theta_1)+\mu_3f(\theta_1+\theta_2)+...+\mu_nf(\theta_1+\theta_2+...+\theta_{n-1})=0, \nonumber \\
&&\mu_3f(\theta_2)+\mu_4f(\theta_2+\theta_3)+...+\mu_1f(\theta_2+\theta_3+...+\theta_{n})=0, \nonumber \\
&&\mu_4f(\theta_3)+\mu_5f(\theta_3+\theta_4)+...+\mu_2f(\theta_3+\theta_4+...+\theta_{n}+\theta_1)=0,\nonumber\\
&&...\label{sistemageral}\\
&&\mu_nf(\theta_{n-1})+\mu_1f(\theta_{n-1}+\theta_{n})+...+\mu_{n-2}f(\theta_{n-1}+\theta_n+\theta_1+...+\theta_{n-3})=0, \nonumber\\
&&\mu_1f(\theta_{n})+\mu_2f(\theta_{n}+\theta_{1})+...+\mu_{n-1}f(\theta_{n}+\theta_1+\theta_2+...+\theta_{n-2})=0, \nonumber\\
&&\theta_1+...+\theta_n=2\pi,\nonumber
\end{eqnarray}
where $f(x) = \ds\sin(x)\left(1 - \frac{1}{8|\sin^3(x/2)|} \right).$
\\
\begin{definition}
We say that a solution $(\tet1,...,\tet2)$ of the system (\ref{sistemageral}) is a central configuration of the planar $1+n$ body problem associated to the masses $\mu_1,...,\mu_n$.
\end{definition} The following results can be found in \cite{albouy} and exhibit the properties of the function $f$ used to prove our results.

\begin{lemma}\label{lema1}
 The function
 $$f(x)=\sin(x)\left(1 - \frac{1}{8|\sin^3(x/2)|} \right), \ \ x\in (0,2\pi)$$
satisfies:
\begin{enumerate}[i)]
 \item $f(\pi-x)=-f(\pi+x), \forall x \in (0,\pi);$
 \item $\displaystyle f'(x)=\cos(x)+\frac{3+\cos(x)}{16|\sin^3(x/2)|}\geq f'(\pi)=-7/8$, for all $x\in(0,2\pi);$
 \item $f'''(x)> 0$, for all $x\in(0,2\pi);$
\item In $(0,\pi)$ there is a unique critical point $\tet c$ of $f$ such that $\tet c>3\pi/5$,  $f'(\te)>0$ in $(0,\tet c)$ and $f'(\te)<0$ in $(\tet c,\pi).$
\end{enumerate}
\end{lemma}

\begin{lemma}\label{lema2albouy}
Consider four points $t_1^L,t_1^R,t_2^L,t_2^R$ such that $0<t_1^L<t_2^L<\te_c<t_2^R<t_1^R<2\pi, f(t_1^L)=f(t_1^R)=f_1$ and $f(t_2^L)=f(t_2^R)=f_2$. Then
$t_2^L+t_2^R<t_1^L+t_1^R.$
\end{lemma}

\begin{cor}\label{corolarioalbouy}
Consider $0<t_1<\theta_c<t_2<2\pi.$ If $f(t_1)\geq f(t_2)$ then $t_1+t_2>2\theta_c>6\pi/5.$
\end{cor}

\begin{figure}
 \centering
  \includegraphics[width=10cm]{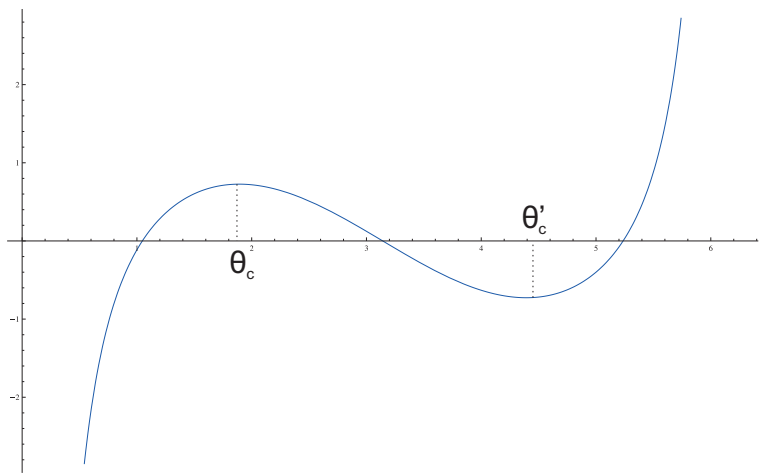}
  \caption{Graphic of $f$.}
\end{figure}

\section{Adding one satellite }\label{addone}

\begin{figure}
  \centering
  \includegraphics[width=8cm]{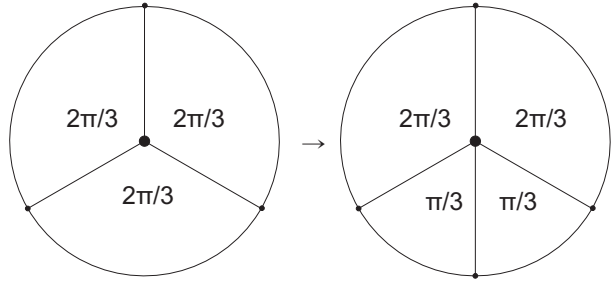}
  \caption{The only stacked central configuration according to Theorema \ref{teoremastacked1}.}
  \label{fig_ccp34}
\end{figure}

\hspace{0,5cm}In this section we consider the planar problem of $1+n$ bodies, where the satellites do not necessarily have the same mass. We are interested in the stacked central configurations when one new satellite is added. In the following result we show that such construction is possible only in one case. See Fig. \ref{fig_ccp34}.

\begin{theorem}\label{teoremastacked1}
Let $(\theta_1,...,\theta_{n-1},\theta_n)$  be a planar central configuration of the $1+n$ body problem where the $n$ satellites have masses $\mu_i, i=1,...,n.$ Suppose that we add another satellite with mass $\mu_{n+1}$ forming a new central configuration of the $1+N$ body problem, with $N=n+1$.
Then $n=3, \theta_1=\theta_2=\theta_3=2\pi/3$, the fourth satellite is on the bisector of one of these angles and the masses of the satellites are all the same.
\end{theorem}

\dem First we will treat the case $n=2$. By hypothesis $(\theta_1,\theta_2)$ is a central configuration relative to the masses
$\mu_1$ and $\mu_2.$ Without lost of generality, suppose that the third satellite, with mass $\mu_3$, is placed such that $(\theta_1,\theta_2',\theta_3')$ is a central configuration of the $1+3$ body problem, where $\theta_2=\theta_2'+\theta_3'$. So the following equations are satisfied:
$$\mu_2f(\theta_1)=0$$
$$\mu_1f(\theta_2)=0$$
$$\mu_2f(\theta_1)+\mu_3f(\theta_1+\theta_2')=0$$
$$\mu_3f(\theta_2')+\mu_1f(\theta_2'+\theta_3')=0$$
$$\mu_1f(\theta_3')+\mu_2f(\theta_3'+\theta_1)=0$$

It is easy to see that the above equations give us
$$f(\theta_1)=f(\theta_1+\theta_2')= f(\theta_2')=0.$$

The roots of $f$ in $(0,2\pi)$ are $\pi/3,\pi$ and $5\pi/3$, which makes impossible the equalities above.

Consider now the case $n=3.$ In the same way, suppose that the fourth satellite, with mass $\mu_4$, is placed between the first satellite and the third one, forming the central configuration $(\theta_1, \theta_2, \theta_3', \theta_4'),$ with $\theta_3=\theta_3'+\theta_4'.$ The following equations are satisfied:
$$\mu_2f(\theta_1)+\mu_3f(\theta_1+\theta_2)=0,$$
$$\mu_3f(\theta_2)+\mu_1f(\theta_2+\theta_3)=0,$$
$$\mu_1f(\theta_3)+\mu_2f(\theta_1+\theta_3)=0,$$
$$\mu_2f(\theta_1)+\mu_3f(\theta_1+\theta_2)+\mu_4f(\theta_1+\theta_2+\theta_3')=0,$$
$$\mu_3f(\theta_2)+\mu_4f(\theta_2+\theta_3')+\mu_1f(\theta_2+\theta_3'+\theta_4')=0,$$
$$\mu_4f(\theta_3')+\mu_1f(\theta_3'+\theta_4')+\mu_2f(\theta_3'+\theta_4'+\theta_1)=0.$$

Since $\theta_3=\theta_3'+\theta_4'$ the equations above give us
$$f(\theta_3')=f(\theta_3'+\theta_2)=f(\theta_3'+\theta_2+\theta_1)=0.$$

So $\theta_3'=\pi/3, \theta_3'+\theta_2=\pi$ and $\theta_3'+\theta_2+\theta_1=5\pi/3.$ Hence, $\theta_1=\theta_2=\theta_3=2\pi/3$ and $\theta_3'=\theta_4'=\pi/3$. Substituting these values into the above equations, we see that $\mu_1=\mu_2=\mu_3=\mu_4$.

Finally we treat the case $n\geq 4.$ Assume that $(\theta_1,...,\theta_{n-1},\theta_n)$ and $(\theta_1,...,\theta_{n-1},\theta_n',\theta_{n+1}'),$ with $\theta_n=\theta_n'+\theta_{n+1}',$ are central configurations associated to the masses $\mu_1,...,\mu_n$ and $\mu_1,...,\mu_n,\mu_{n+1}$, respectively. In fact we are assuming, without lost of generality, that the new satellite, with mass $\mu_{n+1}$, was placed between the first satellite and the $n$th one.

The first equations of the respective systems are
$$\mu_2f(\theta_1)+\mu_3f(\theta_1+\theta_2)+...+\mu_nf(\theta_1+...+\theta_{n-1})=0$$
and
$$\mu_2f(\theta_1)+\mu_3f(\theta_1+\theta_2)+...+\mu_nf(\theta_1+...+\theta_{n-1})+\mu_{n+1}f(\theta_1+...+\theta_{n-1}+\theta_n')=0.$$
This gives us
$$f(\theta_1+...+\theta_{n-1}+\theta_n')=-f(\theta_{n+1}')=0.$$

In the same way, the second equations are
$$\mu_3f(\theta_2)+...+\mu_nf(\theta_2+...+\theta_{n-1})+\mu_1f(\theta_2+...+\theta_{n-1}+\theta_n)=0$$
and
$$\mu_3f(\theta_2)+...+\mu_nf(\theta_2+...+\theta_{n-1})+\mu_{n+1}f(\theta_2+...+\theta_{n-1}+\theta_n')+\mu_1f(\theta_2+...+\theta_{n-1}+\theta_n'+\theta_{n+1}')=0.$$
So
$$f(\theta_2+...+\theta_{n-1}+\theta_n')=-f(\theta_{n+1}'+\theta_1)=0.$$

Analogously, comparing the $i$th equations of the respective systems, we have
$$f(\theta_{n+1}'+\theta_1+\theta_2+...+\theta_{i-1})=0.$$

Thus, we found $n$ distinct roots of $f$ in $(0,2\pi)$ given by
$$\theta_{n+1}',\theta_{n+1}'+\theta_1,\theta_{n+1}'+\theta_1+\theta_2,...,\theta_{n+1}'+\theta_1+\theta_2+...+\theta_{n-1}.$$

So, we have a contradiction for $n\geq 4$ and the Theorem follows.

\cqd
\section{Adding two satellites}\label{addtwo}

\begin{lemma}\label{lema_stacked2} Let $C=\{x_1,...,x_m\}$, with $x_i>0$ and $x_1+...+x_m\leq 2\pi$. Consider four non-empty subsets of $C$, $A_1,A_2,B_1,B_2$ such that $A_1\cup A_2=B_1\cup B_2=C$ and $B_1\subsetneq A_1$. Suppose that
$$f\left(\sum_{x_i\in A_1}x_i\right)=f\left(\sum_{x_i\in A_2}x_i\right)$$
and
$$f\left(\sum_{x_i\in B_1}x_i\right)=f\left(\sum_{x_i\in B_2}x_i\right).$$

Then
$$\sum_{x_i\in C}x_i>2\theta_c>6\pi/5$$ and exactly one of the situations bellow must occur:
\begin{itemize}
  \item $\displaystyle\sum_{x_i\in A_1}x_i=\sum_{x_i\in A_2}x_i$;
  \item $\displaystyle\sum_{x_i\in B_1}x_i=\sum_{x_i\in B_2}x_i$;
  \item $\displaystyle\sum_{x_i\in A_1}x_i=\sum_{x_i\in B_2}x_i.$
\end{itemize}
\end{lemma}

\dem By hypothesis we can use Lemma \ref{lema2albouy} and, observing that
$$\displaystyle\sum_{x_i\in A_1}x_i+\displaystyle\sum_{x_i\in A_2}x_i=\displaystyle\sum_{x_i\in B_1}x_i+\displaystyle\sum_{x_i\in B_2}x_i=\displaystyle\sum_{x_i\in C}x_i,$$
we get
$$\displaystyle\sum_{x_i\in A_1}x_i=\displaystyle\sum_{x_i\in A_2}x_i$$ or $$\displaystyle\sum_{x_i\in B_1}x_i=\displaystyle\sum_{x_i\in B_2}x_i.$$
Since $B_1\subsetneq A_1$ and hence $A_2\subsetneq B_2$, only one of the equalities above can be satisfied. So, by Corollary \ref{corolarioalbouy},
$$\sum_{x_i\in C}x_i>2\theta_c>6\pi/5.$$
\cqd

Considering all satellites with the same mass, the three following propositions give us all possible cases of getting a new planar central configuration of $1+(n+2)$ bodies by adding two satellites to a planar central configuration of $1+n$ bodies.

\begin{figure}
  \centering
  \includegraphics[width=8cm]{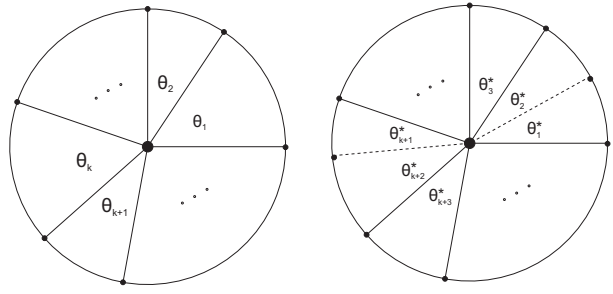}
  \caption{Illustration of the situation treated in Proposition \ref{prop_stacked1}}
  \label{fig_stacked1}
\end{figure}

\begin{prop}\label{prop_stacked1}
Let $(\te_1,...,\te_n)$ be a planar central configuration of $1+n$ bodies, with $n\geq4$ and satellites with the same mass. Suppose that we add two new satellites separated on the circle by at least two of the original satellites. Then the new  $1+(n+2)$ bodies cannot form a central configuration.
\end{prop}
\dem
Without lost of generality, suppose that the new satellites are placed one between the first and second satellites and the other between the $k$th and $(k+1)$th satellites, where $3\leq k\leq n-1$. See Fig. \ref{fig_stacked1}. So,  we get the central configurations
$$(\theta_1,\theta_2,...,\theta_k,...,\theta_n)$$ and
$$(\theta_1^*,\theta_2^*,...,\theta_k^*,\theta_{k+1}^*,\theta_{k+2}^*,...,\theta_{n+2}^*)$$
of $1+n$ and $1+(n+2)$ bodies, respectively, satisfying
$$\theta_1=\theta_1^*+\theta_2^*,$$
$$\theta_k=\theta_{k+1}^*+\theta_{k+2}^*,$$
$$\theta_i=\theta_{i+1}^* \mbox{ for } 2\leq i\leq k-1,$$
and
$$\theta_i=\theta_{i+2}^* \mbox{ for } k+1\leq i \leq n.$$

The respective systems are
\begin{eqnarray}
&&f(\theta_1)+f(\theta_1+\theta_2)+...+f(\theta_1+\theta_2+...+\theta_{n-1})=0, \nonumber \\
&&f(\theta_2)+f(\theta_2+\theta_3)+...+f(\theta_2+\theta_3+...+\theta_{n})=0, \nonumber \\
&&f(\theta_3)+f(\theta_3+\theta_4)+...+f(\theta_3+\theta_4+...+\theta_{n}+\theta_1)=0,\nonumber\\
&&...\label{sistema_staked1}\\
&&f(\theta_{n-1})+f(\theta_{n-1}+\theta_{n})+...+f(\theta_{n-1}+\theta_n+\theta_1+...+\theta_{n-3})=0, \nonumber\\
&&f(\theta_{n})+f(\theta_{n}+\theta_{1})+...+f(\theta_{n}+\theta_1+\theta_2+...+\theta_{n-2})=0 \nonumber
\end{eqnarray}
and
\begin{eqnarray}
&&f(\theta_1^*)+f(\theta_1^*+\theta_2^*)+...+f(\theta_1^*+\theta_2^*+...+\theta_{n+1}^*)=0, \nonumber \\
&&f(\theta_2^*)+f(\theta_2^*+\theta_3^*)+...+f(\theta_2^*+\theta_3^*+...+\theta_{n+2}^*)=0, \nonumber \\
&&f(\theta_3^*)+f(\theta_3^*+\theta_4^*)+...+f(\theta_3^*+\theta_4^*+...+\theta_{n+2}^*+\theta_1^*)=0,\nonumber\\
&&...\label{sistema_staked2}\\
&&f(\theta_{n+1}^*)+f(\theta_{n+1}^*+\theta_{n+2}^*)+...+f(\theta_{n+1}^*+\theta_{n+2}^*+\theta_1^*+...+\theta_{n-1}^*)=0, \nonumber\\
&&f(\theta_{n+2}^*)+f(\theta_{n+2}^*+\theta_{1}^*)+...+f(\theta_{n+2}^*+\theta_1^*+\theta_2^*+...+\theta_{n}^*)=0 \nonumber
\end{eqnarray}

The third equation of (\ref{sistema_staked2}) together with the second equation of (\ref{sistema_staked1}) give us
$$f(\te_2+...+\te_{k-1}+\te_{k+1}^*)+f(\te_2+...+\te_n+\te_1^*)=0\implica$$
$$f(\te_2+...+\te_{k+1}^*)=f(\te_2^*).$$

The $(k+1)$th equation of (\ref{sistema_staked2}) together with the $k$th of (\ref{sistema_staked1}) give us
$$f(\te_{k+1}^*)+f(\te_k+...+\te_n+\te_1^*)=0\implica$$ $$f(\te_{k+1}^*)=f(\te_2^*+\te_2+...+\te_{k-1}).$$

So, by Lemma \ref{lema_stacked2}, we have
$$\te_2^*+\te_2+...+\te_{k-1}+\te_{k+1}^*>\pi.$$

The $(k+3)$th equation of (\ref{sistema_staked2}) together with the $(k+1)$th equation of (\ref{sistema_staked1}) and the first equations of (\ref{sistema_staked1}) and (\ref{sistema_staked2}) give us, respectively
$$f(\te_{k+1}+...+\te_n+\te_1^*)+f(\te_{k+1}+...+\te_n+\te_1+...+\te_{k-1}+\te_{k+1}^*)=0$$
and
$$f(\te_1^*)+f(\te_1^*+\te_2+...+\te_{k-1}+\te_{k+1}^*)=0.$$
So
$$ f(\te_{k+1}+...+\te_n+\te_1^*)=f(\te_{k+2}^*)$$
and
$$f(\te_1^*)=f(\te_{k+2}^*+\te_{k+1}+...+\te_n).$$

Again, by Lemma \ref{lema_stacked2}
$$\te_1^*+\te_{k+2}^*+\te_{k+1}+...+\te_n>\pi.$$

Therefore we have the contradiction
$$\te_1+...+\te_n>2\pi.$$
\cqd

\begin{figure}
  \centering
  \includegraphics[width=8cm]{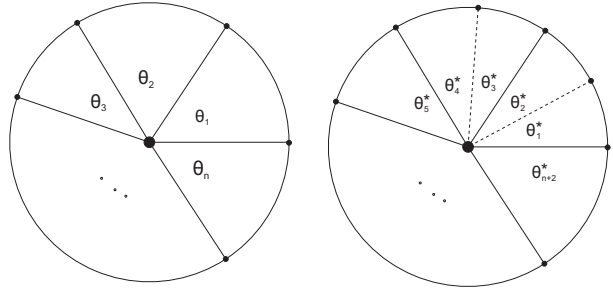}
  \caption{Illustration of the situation treated in Proposition \ref{prop_stacked2}}
  \label{fig_stacked2}
\end{figure}

\begin{figure}
  \centering
  \includegraphics[width=8cm]{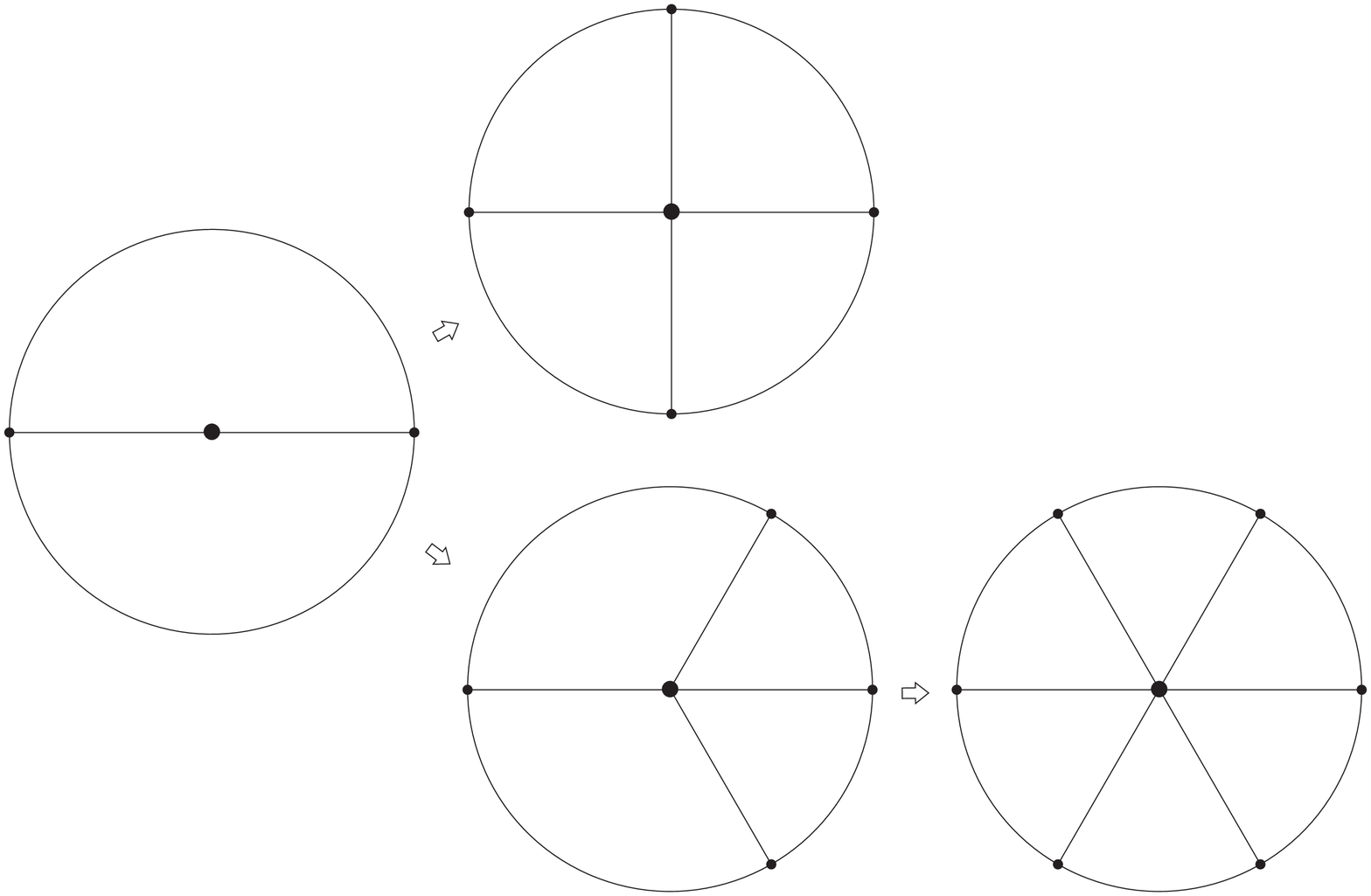}
  \caption{Stackeds central configurations according to Proposition \ref{prop_stacked2}}
  \label{fig2_stacked2}
\end{figure}

\begin{prop}\label{prop_stacked2}
Let $(\te_1,...,\te_n)$ be a planar central configuration of $1+n$ bodies, with $n\geq2$ and satellites with the same mass. Suppose that we add two satellites separated by one of the original satellites. Then the new configuration with $1+(n+2)$ bodies is a central configuration exactly in the cases bellow. See Fig. \ref{fig2_stacked2}.

\noindent $n=2:$\\
Original configuration: collinear = $(\pi,\pi)$.\\
New configurations:
$$\mbox{square = }\left(\frac{\pi}{2},\frac{\pi}{2},\frac{\pi}{2},\frac{\pi}{2}\right), $$
or
$$\mbox{kite = }\left(\frac{\pi}{3},\frac{\pi}{3},\frac{2\pi}{3},\frac{2\pi}{3}\right).$$
$n=4:$\\
Original configuration: kite = $\left(\frac{\pi}{3},\frac{\pi}{3},\frac{2\pi}{3},\frac{2\pi}{3}\right).$\\
New configuration: regular 6-gon = $\left(\frac{\pi}{3},\frac{\pi}{3},\frac{\pi}{3},\frac{\pi}{3},\frac{\pi}{3},\frac{\pi}{3}\right)$.
\end{prop}

\dem First, consider the case $n=2.$ Suppose that we have the central configuration $(\tet1,\tet2)$ and $(\tet1^*,\tet2^*,\tet3^*,\tet4^*)$ of $1+2$ and $1+4$ bodies, respectively, where $\tet1=\tet1^*+\tet2^*$ and $\tet2=\tet3^*+\tet4^*$. The equations are:
\begin{eqnarray}
f(\tet1^*)&=&f(\tet4^*),\label{eq1n2n4}\\
f(\tet3^*)&=&f(\tet2^*),\label{eq2n2n4}\\
f(\tet2^*)+f(\tet2^*+\tet3^*)&=&f(\tet1^*).\label{eqn2n4}
\end{eqnarray}

By Corollary \ref{corolarioalbouy} the equations (\ref{eq1n2n4}) and (\ref{eq2n2n4}), give us
$$\tet1^*=\tet4^* \mbox{ or } \tet3^*=\tet2^*.$$

In fact, by Proposition 2 in \cite{albouy} we have
$$\tet1^*=\tet4^* \mbox{ and } \tet3^*=\tet2^*.$$
So $\tet1=\tet2=\pi.$

The equation (\ref{eqn2n4}) becomes $h(\tet2^*)=0,$ where
$$h(x):=f(x)+f(2x)-f(\pi+x), \ x\in (0,\pi).$$
It is easy to see that $h'''(x)>0$, because $f'''(x)>f'''(\pi+x)>0$. Hence, $h$ has at most three roots in $(0,\pi).$ We can check that
$h(\pi/2)=h(\pi/3)=h(2\pi/3)=0.$
Therefore, $\tet2^*=\pi/2,$ $\tet2^*=\pi/3$ or $\tet2^*=2\pi/3.$ The first root gives us the square configuration and the last two roots correspond to the same configuration, the kite one, given by $(\pi/3,\pi/3,2\pi/3,2\pi/3).$

Consider now $n\geq 3$ and suppose that the two new satellites were placed separated by the original second satellite, as in Fig. \ref{fig_stacked2}. Suppose that we get the central configurations $(\te_1,...,\te_n)$ and $(\te_1^*,...,\te_{n+2}^*)$ satisfying
$$\te_1=\te_1^*+\te_2^*,$$
$$\te_2=\te_3^*+\te_4^*,$$
$$\te_i=\te_{i+2}^*, \mbox{ for } 3\leq i \leq n.$$

Let us see the case $n=3$. The equations reduce to
\begin{eqnarray}
   f(\te_1^*)&=&f(\te_3+\te_4^*),\label{eqprop2st1} \\
	f(\te_3+\te_1^*)&=&f(\te_4^*),\label{eqprop2st2}\\
	f(\te_2^*)&=&f(\te_3^*),\label{eqprop2st3}\\
f(\tet2^*)+f(\tet2^*+\tet3^*)&=&f(\tet3+\tet1^*)+f(\tet1^*).\label{eqprop2st4}
\end{eqnarray}

By Lemma \ref{lema_stacked2}, (\ref{eqprop2st1}) and (\ref{eqprop2st2}) we have
$$\te_1^*+\te_4^*+\te_3>2\te_c>6\pi/5.$$

It follows from this fact and from (\ref{eqprop2st3}) that $\te_2^*=\te_3^*$.

By Lemma \ref{lema_stacked2}, from (\ref{eqprop2st1}) and (\ref{eqprop2st2}) we have exactly one of the following situations
\begin{eqnarray}
\nonumber \te_1^*=\te_4^*,\\
\te_1^*=\te_3+\te_4^*,\\
\nonumber \te_4^*=\te_3+\te_1^*.
\end{eqnarray}
Suppose $\te_1^*=\tet3+\tet4^*$. The equation (\ref{eqprop2st2}) becomes
\begin{equation}\label{eqprop2st5}
f(\tet4^*)=f(\tet4^*+2\tet3).
\end{equation}

Since $(\tet1,\tet2,\tet3)$ is a planar central configuration of the $1+3$ body problem, we know that $\tet1=\tet2, \tet1=\tet3$ or $\tet2=\tet3.$ See \cite{llibre2}. As $\tet1=\tet1^*+\tet2^*=\tet3+\tet4^*+\tet2^*>\tet3,$ we have
$$\tet1=\tet2>\tet3 \mbox{ or } \tet1>\tet2=\tet3.$$
In the first case, we have $\tet3>4\pi/9 \mbox{ and } \tet1<7\pi/9.$ From (\ref{eqprop2st5}) we have $\pi/4<\tet4^*<\pi/3$ and so by Lemma \ref{lema2albouy}, (\ref{eqprop2st5}) implies
$$2\tet4^*+2\tet3>\pi/3+\pi \implica \tet1^*=\tet4^*+\tet3>2\pi/3\implica \tet2^*=\tet1-\tet1^*<7\pi/9-2\pi/3=\pi/9.$$
So $f(\tet2^*)<f(\pi/9)<-7.$ Therefore (\ref{eqprop2st4}) is not satisfied.

If $\tet1>\tet2=\tet3$ then $\tet3<5\pi/18.$ Equation (\ref{eqprop2st5}) gives us $\tet4^*>\pi/3$ and hence $\tet3=\tet2=\tet3^*+\tet4^*>\pi/3>5\pi/18>\tet3.$

Analogously we see that the case $\tet4^*=\tet3+\tet1^*$ is impossible too. So $$\te_1^*=\te_4^*$$
and, consequently
$$\te_1=\te_2.$$
So, (\ref{eqprop2st1}) becomes
\begin{equation}\label{eqp2st}
f(\tet1^*)=f(\tet3+\tet1^*).
\end{equation}

Since $(\te_1,\te_1,\te_3)$ is a central configuration of $1+3$ bodies, we see that
$$\f1+f(2\tet1)=0$$
and we will consider the three roots of this equation again. See \cite{llibre2}.

If $\tet1=2\pi/3$, then $\tet3 =2\pi/3$ and, by (\ref{eqp2st}) we have
$$\tet1^*=\pi/3.$$

But, it is easy to see that the resultant configuration does not agree with (\ref{eqprop2st4}).

Consider now that $\tet1$ is the smallest solution of $f(\tet1)+f(2\tet1)=0$. We know that, in this case, $\te_1<50^{\circ}<\tet c/2$, so $\tet3>\tet c'=2\pi-\tet c$, the second critical point of $f$. Therefore, $\tet3$ lies on the last increasing interval of $f$. Thus,
$$f(\tet1^*)<f(\tet1)=f(\tet3)<f(\tet3+\tet1^*),$$
contradicting (\ref{eqp2st}).

Now, let $\tet1$ be the biggest solution of $f(\tet1)+f(2\tet1)$. We know that $138^{\circ}<\tet1<139^{\circ}$ and thus, $$\te_3>82^{\circ}.$$

It is easy to see that, in this case, by (\ref{eqp2st}), we have $\pi/3<\tet1^*<\tet c<\tet1^*+\tet3<\pi.$

If $\tet1^*\geq 73^{\circ}$ then $\tet1^*+\tet3>155^{\circ}$, what implies
$$0,365...=f(155^{\circ})>f(\tet1^*+\tet3)=f(\tet1^*)>f(73\pi/180)=0,388...$$
Therefore,
$$\tet1^*<73^{\circ}$$
and consequently
$$\tet2^*=\tet1-\tet1^*>65^{\circ}.$$

The function $f(\te)+f(2\te)$ is increasing from $0$ to its first critical point and then decreasing to its second critical point, bigger than $2\pi/3$. Since $f'(65^{\circ})+2f'(130^{\circ})>0$, the first critical point is bigger than $65^{\circ}.$ Thus,
$$f(\tet2^*)+f(2\tet2^*)>f(65^{\circ})+f(130^{\circ})>0,8.$$

However, $2f(\tet1^*)<2f(73^{\circ})<2(0,39)<0,8$
contradicting (\ref{eqprop2st4}).

Now we consider $n\geq5$. The system for the $\tet i^*$'s, by the relations between the angles and the properties of $f$, reduces to
\begin{eqnarray}
 \label{eqsta6} f(\te_4^*) &=& f(\te_1^*+\te_3+...+\te_n), \\
 \label{eqsta7} f(\te_4^*+\te_3) &=&  f(\te_1^*+\te_4+...+\te_n),\\
 \label{eqsta8} f(\te_1^*) &=&  f(\te_4^*+\te_3+...+\te_n),\\
 \label{eqsta9} f(\te_1^*+\te_n) &=&  f(\te_4^*+\te_3+...+\te_{n-1}).
\end{eqnarray}

%
The equations (\ref{eqsta6}) and (\ref{eqsta8}), by Lemma \ref{lema_stacked2}, imply exactly one of following situations:
\begin{eqnarray}\label{eqsta10}
\nonumber  \te_1^* &=& \te_4^*, \\
  \te_1^* &=& \te_4^*+\te_3+...+\te_n,\mbox{ or } \\
\nonumber  \te_4^* &=& \te_3+...+\te_n+\te_1^*.
\end{eqnarray}
Analogously, the equations (\ref{eqsta6}) and (\ref{eqsta7}), by Lemma \ref{lema_stacked2}, imply that exactly one of following equations is verified:
\begin{eqnarray}\label{eqsta11}
\nonumber  \te_4^* &=& \te_1^*+\te_4+...+\te_n, \\
  \te_4^* &=& \te_1^*+\te_3+...+\te_n,\mbox{ or } \\
\nonumber  \te_4^*+\te_3 &=&\te_1^*+ \te_4+...+\te_n
\end{eqnarray}
and, by (\ref{eqsta8}) and (\ref{eqsta7}), the same Lemma implies  exactly one of the following equalities:
\begin{eqnarray}\label{eqsta12}
\nonumber  \te_1^* &=& \te_4^*+\te_3, \\
  \te_1^* &=& \te_4^*+\te_3+...+\te_n, \mbox{ or } \\
\nonumber  \te_4^*+\te_3 &=& \te_4+...+\te_n+\te_1^*.
\end{eqnarray}

The only possibility for a single equation to be satisfied in each group (\ref{eqsta10}), (\ref{eqsta11}) and (\ref{eqsta12}) is
$$\te_1^*=\te_4^*.$$

From the group of equations (\ref{eqsta12}) we have
\begin{equation}\label{stackedte3}
    \te_3=\te_4+...+\te_n.
\end{equation}

Again, by Lemma \ref{lema_stacked2} applied to (\ref{eqsta8}) and (\ref{eqsta9}) and, since $\te_1^*=\te_4^*$, we have
$$\te_n = \te_3+...+\te_{n-1},$$
contradicting (\ref{stackedte3}).

Finally we consider the case $n=4$. The equations are
\begin{eqnarray}
f(\tet1^*)&=&f(\tet3+\tet4+\tet4^*),\label{stackedn4eq1}\\
f(\tet2^*)+f(\tet2^*+\tet3^*)+f(\tet2^*+\tet2)&=&f(\tet4+\tet1^*)+f(\tet1^*),\label{stackedn4eq2}\\
f(\tet3^*)&=&f(\tet2^*),\label{stackedn4eq3}\\
f(\tet4^*)&=&f(\tet3+\tet4+\tet1^*),\\
f(\tet4+\tet1^*)&=&f(\tet4^*+\tet3).
\end{eqnarray}

Using the same argument of the case $n\geq 5$ we get
$$\tet1^*=\tet4^*$$ and $$\tet3=\tet4.$$

Moreover, from (\ref{stackedn4eq1}) and (\ref{stackedn4eq3}), for the angles no greater than $2\pi$ we must have
$$\tet2^*=\tet3^*.$$

In this case, the original configuration with $n=4$, has a symmetry axis passing through two satellites and hence, it should be the square $\tet1=\tet2=\tet3=\tet4=\pi/2$, or the kite configuration given by
$$\tet1=\tet2=\pi/3 \mbox{ and } \tet3=\tet4=2\pi/3$$ or
$$\tet1=\tet2=2\pi/3 \mbox{ and } \tet3=\tet4=\pi/3.$$

If $\tet3=2\pi/3$, from (\ref{stackedn4eq1}), we have
$$f(\tet1^*)=f(\tet1^*+4\pi/3).$$

So,
$$\tet1^*>\pi/4\implica \tet2^*=\tet1-\tet1^*=\pi/3-\tet1^*<\pi/12$$
$$ \implica f(\tet2^*)<f(\pi/12)<-14.$$

It is easy to see that the equation (\ref{stackedn4eq2}) is not satisfied in this case.

If $\tet3=\pi/3$, from (\ref{stackedn4eq1}) we have
$$f(\tet1^*)=f(\tet1^*+2\pi/3).$$

Therefore, by the graph of $f$
$$\tet1^*=\pi/3$$
$$\implica \tet2^*=\pi/3$$
hence we have a regular 6-gon.

Finally the case $\tet1=\tet2=\tet3=\tet4=\pi/2.$ Taking $x=\tet1^*$, we can observe that the equations are
\begin{eqnarray}
f(x)&=&f(\pi+x),\label{quadradoparan6eq1}\\
3f(x)&=&f(\pi/2-x)+f(\pi-2x)\label{quadradoparan6eq2}.
\end{eqnarray}

We will show that the roots of the equations above are different. We can write
$$f(x)=\frac{1}{4}\cos(x/2)\sen^{-2}(x/2)(8\sen^3(x/2)-1).$$

It is easy to see that, from (\ref{quadradoparan6eq1}), we have $0<x<\pi/3.$
Taking $u=\sen(x/2)$ and $v=\cos(x/2)$ we get
\begin{eqnarray*}
4f(x)&=&vu^{-2}(8u^3-1),\\
4f(\pi+x)&=&-uv^{-2}(8v^3-1),\\
4f(\pi-2x)&=&-2uv(v^2-u^2)^{-2}(8(v^2-u^2)^3-1),\\
4f(\pi/2-x)&=&\sqrt{2}(u+v)(v-u)^{-2}(2\sqrt{2}(v-u)^3-1).
\end{eqnarray*}

Setting $P=4f(x)-4f(\pi+x)$ and $Q=12f(x)-4f(\pi/2-x)-4f(\pi-2x),$ the equations are
$$P=0, \ \ Q=0.$$

Denote by $p$ and $q$ the numerators of $P$ and $Q$, respectively. Then
$$p(u,v)=-u^3 - v^3 + 16 u^3 v^3$$ and
$$q(u,v)=\sqrt{2} u^5+4 u^8-2 u^3 v-3 u^4 v+3 \sqrt{2} u^4 v+24 u^7 v-16 u^9 v+3 \sqrt{2} u^3 v^2-12 u^6 v^2+6 u^2 v^3+$$
$$+\sqrt{2} u^2 v^3-48 u^5 v^3+48 u^7 v^3+12 u^4 v^4-3 v^5+24 u^3 v^5-48 u^5 v^5-4 u^2 v^6+16 u^3 v^7.$$

Let $R_1$ and $R_2$ be the resultants of $p$ and $q$ with respect to $u$ and $v$, respectively. $R_1$ and $R_2$ are polynomials in the variables $v$ and $u$, respectively, such that, if $(u_0,v_0)$ is a solution for $p=0$ and $q=0$, then $R_1(v_0)=R_2(u_0)=0.$

We get
$$R_1(v)=-2v^{12}r_1(v),$$
where
\begin{eqnarray*}
r_1(v)&=&-4+18 v+18 \sqrt{2} v-135 v^2+384 v^3-1800 v^4-1728 \sqrt{2} v^4+ 12240 v^5+...\\&&...+231928233984 v^{32}-21474836480 v^{33}-154618822656 v^{34}+34359738368 v^{36}
\end{eqnarray*}
and $$R_2(u)=2u^{12}r_2(u),$$
where
\begin{eqnarray*}
r_2(u)&=&4+18 u-18 \sqrt{2} u+135 u^2-384 u^3-1800 u^4+1728 \sqrt{2} u^4-...\\&&...-21474836480 u^{33}-154618822656 u^{34}+34359738368 u^{36}.
\end{eqnarray*}

We reject the solution $u=v=0$ and we focus on polinomials $r_1(v)$ and $r_2(u).$ We will show that $r_1(v)$ have no roots in $[9/10,1]$. The substitution $v=\frac{9+s}{10+s}$  maps $(0,\infty)$ into $(9/10,1)$ and we get
\begin{eqnarray*}
r_1(s)&=&\left.\frac{1}{(10+s)^{36}}\right(-1008258045536460519041298702036036780097011712+\\
&&+703154989278597577665542304052850511052800000 \sqrt{2}-...\\
&&...-6975267000684 s^{35}+4183796399670 \sqrt{2} s^{35}-21095906033 s^{36}+\\&&+\left.12595023450 \sqrt{2} s^{36}\right).
\end{eqnarray*}
$r_1(s)$ has no roots in $(0,\infty)$ since its numerator is a polynomial of degree 36 where all terms are negative. So we have $0<v<9/10.$


As $u^2+v^2=1,$ and $u=\sin(x/2)\leq 1/2$ then
$$\frac{21}{50}<\frac{\sqrt{19}}{10}=\sqrt{1-(9/10)^2}<u\leq \frac{1}{2}.$$

Taking the rational substitution $u=\frac{21+t}{50+2t}$, we have


\begin{eqnarray*}
r_2(t)&=&\left.\frac{1}{16(25+t)^{36}}\right(39355704304464511525941168410542220606454431557393-\\
&&-39249115788439924325393189378990445455012929687500 \sqrt{2}-...\\
&&...\left.-57616560 \sqrt{2} t^{34}+168012 t^{35}-152568 \sqrt{2} t^{35}+217 t^{36}-196 \sqrt{2} t^{36}\right).
\end{eqnarray*}

Again, all coefficients of the numerator are negative. This shows that $r_2(u)$ does not have roots in $\left(\frac{21}{50},\frac{1}{2}\right)$. Therefore the solutions of (\ref{quadradoparan6eq1}) and (\ref{quadradoparan6eq2}) are different and the proposition follows.
\cqd

\begin{figure}
  \centering
  \includegraphics[width=8cm]{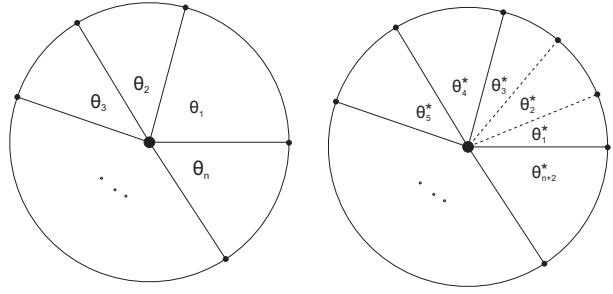}
  \caption{Illustration of the situation treated in Proposition \ref{prop_stacked3}}
  \label{fig_stacked3}
\end{figure}

\begin{figure}
  \centering
  \includegraphics[width=8cm]{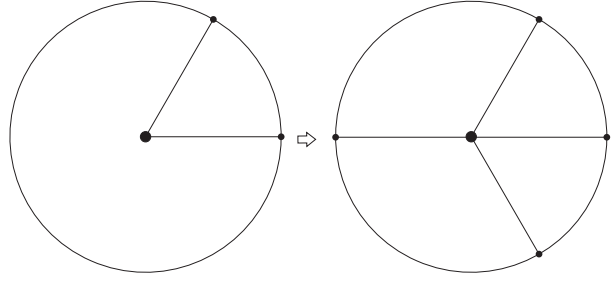}
  \caption{Stacked central configuration according to Proposition \ref{prop_stacked3}}
  \label{fig3_stacked2}
\end{figure}

\begin{prop}\label{prop_stacked3}
Let $(\te_1,...,\te_n)$ be a planar central configuration of $1+n$ bodies. Suppose that we add two satellites between two of the consecutive original  satellites and we get a new planar central configuration (with 1+(n+2) bodies). Then $n=2,$ in the original configuration we have an equilateral triangle formed by the satellites and the massive body and the new configuration is a kite given by $\left(\frac{\pi}{3},\frac{\pi}{3},\frac{2\pi}{3},\frac{2\pi}{3}\right).$ See Fig. \ref{fig3_stacked2}.

\end{prop}

\dem
We suppose that the new satellites are placed between the first and the second original satellites.
The relations between the angles are

$$\te_1=\te_1^*+\te_2^*+\te_3^*$$
and
$$\te_i=\te_{i+2}^*, \ 2\leq i \leq n.$$

See Fig. \ref{fig_stacked3}

Firstly we consider the case $n=2.$ The following equations are satisfied:
\begin{eqnarray}
f(\tet1^*)&=&f(\tet2+\tet3^*),\label{eqst1n2n4}\\
f(\tet1^*)&=&f(\tet2^*)+f(\tet2^*+\tet3^*),\label{eqst2n2n4}\\
f(\tet2^*)&=&f(\tet1^*)+f(\tet3^*),\label{eqst3n2n4}\\
f(\tet3^*)&=&f(\tet1^*+\tet2).\label{eqst4n2n4}
\end{eqnarray}

As $(\tet1,\tet2)$ is a planar central configuration for $n=2$ we get $\tet2=\pi,\tet2=5\pi/3$ or $\tet2=\pi/3.$

If $\tet2=\pi,$ then $\tet1=\tet1^*+\tet2^*+\tet3^*=\pi$ and, by Lemma \ref{lema_stacked2}, the equations (\ref{eqst1n2n4}) and (\ref{eqst4n2n4}) imply $\tet1^*=\tet3^*.$
So $f(\tet1^*)=f(\pi+\tet1^*)$ and by the graph of $f$, $\tet1^*<\pi/3$ and, thus $f(\tet1)<0.$

From (\ref{eqst3n2n4}) we get
$$f(\tet2^*)=2f(\tet1^*)<0\implica \tet2^*<\frac{\pi}{3}.$$ This implies the contradiction
$$ \tet1^*+\tet2^*+\tet3^*<\pi=\tet1.$$

If $\tet2=5\pi/3$ we conclude by the same argument that $\tet1^*=\tet3^*$ and so
$$f(\tet1^*)=f(\tet1^*+5\pi/3).$$
By the graph of $f$ this equations is impossible because we should have $\tet1^*>\pi/3$ and the sum of the angles should be bigger than $2\pi.$

Finally the case $\tet2=\pi/3.$ By (\ref{eqst1n2n4}) and (\ref{eqst4n2n4}), the Lemma \ref{lema_stacked2} gives us three disjoint possibilities, namely, $\tet1^*=\pi/3+\tet3^*, \tet3^*=\pi/3+\tet1^*$ or $\tet1^*=\tet3^*.$

If $\tet1=\pi/3+\tet3^*$ then, by (\ref{eqst4n2n4}), we have
$$f(\tet3^*)=f(2\pi/3+\tet3^*)\implica \tet3^*=\pi/3\implica f(\tet3^*)=0.$$
From (\ref{eqst2n2n4}) and (\ref{eqst3n2n4}) we get $f(\tet2^*+\pi/3)=0,$ so $\tet2^*=2\pi/3$ and consequently we have the kite configuration.

The case $\tet3=\pi/3+\tet1^*$ is analogous and results in the same configuration.

Now consider $\tet1^*=\tet3^*$. So $f(\tet1^*)=f(\pi/3+\tet1^*)$ and thus $$\pi/3<\tet1^*<\tet c<\tet1^*+\pi/3<\pi.$$

By Lemma \ref{lema2albouy} we have $$2\tet1^*+\pi/3<\pi/3+\pi\implica \tet1^*<\pi/2$$
and
$$2\tet1^*+\pi/3>2\tet c>6\pi/5 \implica \tet1^*>\frac{13\pi}{30}.$$

From (\ref{eqst3n2n4}) we have $f(\tet2^*)=2f(\tet1^*).$ Because $\tet1^*<\pi/2$, then $\tet2^*=5\pi/3-2\tet1^*>2\pi/3.$ Therefore
$$2f(\tet1^*)=f(\tet2^*)<f(2\pi/3)<0,7.$$
However, one checks that this fact contradicts the inequality $\tet1^*>13\pi/30,$ because $f(13\pi/30)>f(2\pi/5)>0,36.$

Consider the case $n=3.$ The equations become
$$f(\te_1)=f(\te_2)=f(\te_3),$$
$$f(\te_1^*)=f(\te_2+\te_3+\te_3^*),$$
$$f(\te_3^*)=f(\te_2+\te_3+\te_1^*),$$
$$f(\te_3+\te_1^*)=f(\te_2+\te_3^*).$$

By Lema \ref{lema_stacked2}, the three last equations above imply
$$\te_1^*=\te_3^*.$$
Hence
$$\te_2=\te_3.$$

We get then,
\begin{eqnarray}
f(\te_1^*)&=&f(2\te_2+\te_1^*),\label{eqprop3stacked1}\\
f(\te_2)&=&f(2\te_1^*+\te_2^*),\label{eqprop3stacked2}\\
f(\te_2^*)+f(\te_1^*+\te_2^*)&=&f(\te_1^*+\te_2)+f(\te_1^*).\label{eqprop3stacked3}
\end{eqnarray}

On the other hand, the angle $\te_2$ satisfies
$$f(\te_2)+f(2\te_2)=0.$$
This equation  has three solutions in $(0,\pi)$, namely,
$\te_2=2\pi/3, \te_2=a \mbox{ and } \te_2=b,$
with $0<a<50^{\circ} \mbox{ and } 2\pi/3<b<\pi.$ See \cite{llibre2}.

If $\te_2=2\pi/3$ then $\te_1=2\pi/3$. Since $f$ is injective in $(0,\pi/4)$, by (\ref{eqprop3stacked1}) we have $\te_1^* > \pi/4$. So $\te_2^*=\te_1-2\te_1^*=2\pi/3-2\te_1^*<\pi/6.$ Therefore
$$f(\te_2^*)<-3.$$

It is not hard to see that (\ref{eqprop3stacked3}) is not satisfied in this case.

If $\te_2=a<\te_c/2$, then $\te_1>2\pi-\te_c$. Since $a>\pi/4$ then $2\te_2>\pi/2$.

From the graph of $f$ we can conclude that $\pi/3<\te_1^*<\te_c$. Furthermore, we claim that $\te_1^*<7\pi/18$. In fact, if $\te_1^*>7\pi/18$ then $\te_1^*+2\te_2>8\pi/9$. So, we have the contradiction:
$$0.317<f(7\pi/18)<f(\te_1^*)=f(\te_1^*+2\te_2)<f(8\pi/9)<0.298.$$

Hence,
$$\te_1^*<7\pi/18\implica \te_2^*>2\pi/3$$
$$\implica \te_2+\te_1^*<2\pi/3<\te_2^*$$
$$\implica f(\te_2^*)<f(\te_2+\te_1^*).$$

From (\ref{eqprop3stacked3}) we get
$$f(\te_1^*+\te_2^*)>f(\te_1^*)=f(\te_1^*+2\te_2)$$
$$\implica \te_1^*+2\te_2>\te_1^*+\te_2^*$$
$$\implica 2\te_2>\te_2^*,$$

However, this is impossible because $\te_2<50^{\circ}$ and $\te_2^*>120^{\circ}$.

$\te_1=b>138^{\circ}$ is the remaining case. But this assumption leads us to $\te_1<\pi/2$, a contradiction to the inequality $\te_1^*>\pi/4$ which came from (\ref{eqprop3stacked1}) and the injectivity of the function $f$ in $(0,\pi/4)$.

Now suppose $n\geq 4.$ We will show that the new configuration is not a central one.
From
$$f(\te_1)+f(\te_1+\te_2)+...+f(\te_1+...+\te_{n-1})=0,$$
the first equation for the new configuration becomes
\begin{equation}\label{eqtheta1estrela}
f(\te_1^*)=f(\te_3^*+\te_2+...+\te_n).
\end{equation}

Analogously for the other equations, we get:
\begin{equation}\label{eqtheta3estrela}
    f(\te_3^*)=f(\te_2+...+\te_n+\te_1^*),
\end{equation}
\begin{equation}\label{eqtheta3estrela2}
    f(\te_3^*+\te_2)=f(\te_3+...+\te_n+\te_1^*),
\end{equation}
\begin{equation}\label{eqtheta3estrela23}
    f(\te_3^*+\te_2+\te_3)=f(\te_4+...+\te_n+\te_1^*).
\end{equation}

Notice that we can apply Lemma \ref{lema_stacked2} to equations (\ref{eqtheta1estrela}), (\ref{eqtheta3estrela}), (\ref{eqtheta3estrela2}) and (\ref{eqtheta3estrela23}) chosen pairwise. So, from equations (\ref{eqtheta1estrela}) and (\ref{eqtheta3estrela}) we have exactly one of the following equalities:
\begin{eqnarray}\label{eqsta1}
 \nonumber \te_1^* &=& \te_3^*, \\
  \te_1^* &=& \te_3^*+\te_2+...+\te_n,  \\
 \nonumber \te_3^* &=& \te_2+...+\te_n+\te_1^*.
\end{eqnarray}

From the equations (\ref{eqtheta1estrela}) and (\ref{eqtheta3estrela2}) we have exactly one of the following equalities:
 \begin{eqnarray}\label{eqsta2}
  \te_1^* &=& \te_3^*+\te_2, \nonumber\\
  \te_1^* &=& \te_3^*+\te_2+...+\te_n,  \\
  \te_3^*+\te_2 &=& \te_3+...+\te_n+\te_1^*. \nonumber
\end{eqnarray}

And from (\ref{eqtheta3estrela}) and (\ref{eqtheta3estrela2}) we have exactly one of the situations below:
 \begin{eqnarray}\label{eqsta3}
  \te_3^* &=& \te_2+...+\te_n+\te_1^*, \nonumber\\
  \te_3^*+\te_2 &=& \te_1^*+\te_3+...+\te_n, \\
  \te_3^* &=& \te_3+...+\te_n+\te_1^*. \nonumber
\end{eqnarray}

It is easy to see that from (\ref{eqsta1}), (\ref{eqsta2}) and (\ref{eqsta3}), we have
\begin{eqnarray}\label{eqsta4}
\nonumber \te_1^*&=&\te_3^*,\\
\te_2&=&\te_3+...+\te_n.
\end{eqnarray}

Now, as $\te_1^*=\te_3^*$, by Lemma \ref{lema_stacked2} applied to (\ref{eqtheta1estrela}) and (\ref{eqtheta3estrela23}) we have
$$\te_2+\te_3 = \te_4+...+\te_n.$$

However, this contradicts the equation (\ref{eqsta4}) and the Proposition follows.

\cqd

\noindent {\bf Remark:} Note that the argument to show that there exists no stacked central configurations in the situations treated in Propositions \ref{prop_stacked1} and \ref{prop_stacked3} for $n\geq 4$ and for $n\geq 5$ in Proposition \ref{prop_stacked2}, still works in the more general situation where the original satellites do not necessarily have the same masses, but the added ones have equal masses.

\begin{theorem}
Consider a central configuration of the $1+n$ body problem associated to the masses $\mu_1,...,\mu_n$. Suppose that we add two satellites with the same masse $\mu$. If $n\geq 5$ then the new configuration is not a central one. If the two new satellites are not placed separated by exactly one of the original satellites, then the result holds for $n\geq 4.$
\end{theorem}

\dem
The proofs of Propositions \ref{prop_stacked1} and \ref{prop_stacked3} for $n\geq 4$ and Proposition \ref{prop_stacked2} for $n\geq 5$ have used the Lemma \ref{lema_stacked2} in some equations. The equations in the present case are the same as those because the parameter $\mu$ appears now multiplying both sides of the equality, hence it can be canceled. Therefore the proof does not change.
\cqd

\end{document}